\documentclass[12pt, oneside]{article}
\usepackage[utf8]{inputenc}
\usepackage{amsmath}
\usepackage{amsfonts}
\usepackage{amssymb}
\usepackage{amsthm}

\newtheorem{theorem}{Theorem}[section]

\newtheorem{example}[theorem]{Example}

\newtheorem{remark}{\sc Remark}
\newtheorem{lemma}{\sc Lemma}[section]
\newtheorem{corollary}{\sc Corollary}[section]
\newtheorem{definition}{\sc Definition}[section]

\newcommand{\be}{\begin{eqnarray}}
\newcommand{\ee}{\end{eqnarray}}
\newcommand{\Be}{\begin{eqnarray*}}
\newcommand{\Ee}{\end{eqnarray*}}
\newcommand{\bee}{\begin{equation}}
\newcommand{\eee}{\end{equation}}
\newcommand{\ba}{\begin{array}}
\newcommand{\ea}{\end{array}}
\newcommand{\bl}{\begin{lemma}}
\newcommand{\el}{\end{lemma}}
\newcommand{\bd}{\begin{definition}}
\newcommand{\ed}{\end{definition}}
\newcommand{\bt}{\begin{theorem}}
\newcommand{\et}{\end{theorem}}
\newcommand{\bp}{\begin{proof}}
\newcommand{\ep}{\end{proof}}
\newcommand{\bi}{\begin{itemize}}
\newcommand{\ei}{\end{itemize}}
\newcommand{\br}{\begin{remark}}
\newcommand{\er}{\end{remark}}
\newcommand{\bc}{\begin{corollary}}
\newcommand{\ec}{\end{corollary}}
\newcommand{\bex}{\begin{example}}
\newcommand{\eex}{\end{example}}

\begin{document}

\date{}
\title{\textbf{On the Conformal change of a Douglas space of second kind with special $(\alpha, \beta )$-metric}}
\maketitle
\begin{center}
\author{\textbf{Gauree Shanker$^1$, Sruthy Asha Baby$^{2}$ }}
\end{center}

\begin{center}

$^1$Department of Mathematics and Statistics\\[0pt]
School of Basic and Applied sciences,\\[0pt]
Central University of Punjab, Bathinda-151001, India\\[0pt]
Email: gshankar@cup.ac.in \\[0pt]
$^{2}$Department of Mathematics and Statistics\\[0pt]
Banasthali University, Banasthali\\[0pt]
Rajasthan-304022, India\\[0pt]
Email: sruthymuthu123@gmail.com.

\textbf{Abstract}
\end{center}

{\small The notion of a Douglas space of second kind of a Finsler space
with $(\alpha, \beta)$-metric was introduced by I. Y. Lee [\ref{CD-5}].
Since then, so many geometers have studied this topic e. g., [\ref{CD-14}]. In
this paper, we prove that a Douglas space of second kind with special $%
(\alpha, \beta)$-metric $\alpha +\epsilon \beta + k \frac{\beta^2}{\alpha }$
is conformally transformed to a Douglas space of second kind. Further, we
obtain some results which prove that a Douglas space of second kind with
certain $(\alpha, \beta)$-metrics such as Randers metric, Kropina metric,
first approximate Matsumoto metric and Finsler space with square metric is
conformally transformed to a Douglas space of second kind. }\newline



\section{Introduction}

The notion of Douglas space was introduced by S. Bacso and M. Matsumoto [\ref%
{CD-2}] as a generalization of Berwald space from view point of geodesic
equations. Also, they consider the notion of Landsberg space as a
generalization of Berwald space. Recently, the notion of weakly-Berwald
space as another generalization of Berwald space was introduced by S. Bacso
and B. Szilagyi [\ref{CD-11}]. It is remarkable that a Finsler space is a
Douglas space if and only if the Douglas tensor $D^i_{hjk}$ vanishes
identically [\ref{CD-3}].\newline

The theories of Finsler spaces with an $(\alpha, \beta)$-metric and Berwald
spaces with an $(\alpha, \beta)$-metrics ([\ref{CD-16}], [\ref{CD-7}], [\ref%
{CD-15}]) have contributed alot to the development of Finsler geometry [\ref{CD-8}%
]. The conformal theory of Finsler spaces was introduced by M. S. Kneblman [%
\ref{CD-12}] in 1929 and this theory has been investigated in detail by M.
Hashiguchi [\ref{CD-4}]. Later on Y.D. Lee [\ref{CD-6}] and B. N. Prasad [%
\ref{CD-10}] found conformally invariant tensors in the Finsler space with $%
(\alpha, \beta)$-metric under conformal $\beta$-change.\newline

The purpose of the present paper is to prove that a Douglas space of second
kind with a special $(\alpha, \beta)$-metric given by $L = \alpha + \epsilon
\beta + k \frac{\beta^2}{\alpha}$, where $\epsilon$ and k are constants is
conformally transformed to a Douglas space of second kind. Further, we obtain
some results which prove that the Douglas space of second kind with certain $%
(\alpha, \beta)$-metrics such as Randers metric, Kropina metric, first
approximate Matsumoto metric and Finsler space with square metric is
conformally transformed to a Douglas space of second kind.

\section{Preliminaries}

A Finsler space $F^n$ = $(M^n, L(\alpha, \beta))$ is said to be with an $%
(\alpha, \beta)$-metric, if $L(\alpha, \beta)$ is a positively homogenous
function of $\alpha$ and $\beta$ of degree one, where $\alpha^2 =
a_{ij}(x)y^iy^j$ and $\beta = b_i(x)y^i$. The space $R^n = (M^n, \alpha)$ is
called the Riemannian space associated with $F^n$. We shall use the
following symbols [\ref{CD-8}]: 
\begin{eqnarray*}
&&b^i=a^{ir}b_r,\ \ b^2 = a^{rs}b_rb_s \\
&& 2r_{ij} = b_{i:j} + b_{j:i},\ \ 2s_{ij} = b_{i:j} - b_{j:i}, \\
&& s^i_j = a^{ir}s_{rj},\ \ s^i_j = a^{ir}s_{rj},\ \ s_j = b_rs^r_j.
\end{eqnarray*}
The Berwald connection $B\Gamma = \{G^i_{jk}, G^i_j\}$ of $F^n$ plays an
important role in the present paper. Denote by $B^i_{jk}$, the difference
tensor of $G^i_{jk}$ from $\gamma^i_{jk}$:\newline
\begin{equation}
G^i_{jk} (x,y)= \gamma^i_{jk} (x) + B^i_{jk} (x,y).
\end{equation}
With the subscript 0, transvecting by $y^i$, we have 
\begin{equation}
G^i_j = \gamma^i_{oj} + B^i_j \ \ \ and\ \ 2 G^i = \gamma^i_{00} + 2B^i,
\end{equation}
and then $B^i_j = \dot{\partial}_jB^i$ and $B^i_{jk} = \dot{\partial}_k B^i_j
$.

The geodesics of an n-dimensional Finsler space $F^n$ = $(M^n, L)$ are given
by the system of differential equations [\ref{CD-4}] 
\begin{equation}
\frac{d^2x^i}{dt^2}y^j - \frac{d^2x^j}{dt^2}y^i + 2 (G^iy^j - G^jy^i) = 0;
y^i=\frac{dx^i}{dt},
\end{equation}
in a parameter t. The spray function $G^i (x,y)$ is given by 
\begin{equation}
2G^i(x,y)= g^{ij}(y^r\dot{\partial}_j \partial_r F - \partial_j F)
=\gamma^i_{jk}y^jy^k,
\end{equation}
where $\dot{\partial_i} = \frac{\partial }{\partial x^i}$, F = $\frac{L^2}{2}
$, $\gamma^i_{jk}$ are the Christoffel symbols constructed from $g_{ij}(x,y)$
with respect to $x^i$ and $g^{ij}(x,y)$ is the inverse of the metric tensor 
$g_{ij}(x,y)$.

It is well known [\ref{CD-2}] that a Finsler space $F^n$ becomes a Douglas space if and only if
the Douglas tensor 
\begin{equation}
D^h_{ijk} = G^h_{ijk} - \frac{1}{n + 1}(G_{ijk} y^h + G_{ij}\delta^h_k +
G_{jk}\delta^h_i + G_{ki}\delta^h_j),
\end{equation}
vanishes identically, where $G^h_{ijk} = \dot{\partial_k G^h_{ij}}$ is the h$%
\nu$-curvature tensor of the Berwald connection $B\Gamma$.\newline

Further, a Finsler space $F^n$ is said to be a Douglas space [\ref{CD-9}], if \newline
\begin{equation}
D^{ij} = G^i(x,y)y^j - G^j(x,y)y^i,  \label{C.2.1}
\end{equation}
are homogenous polynomials in $(y^i)$ of degree three.

Differentiating (\ref{C.2.1}) by $y^m$ and contracting m and j in
the obtained equation, we have 
\begin{equation}
D^{im}_m = (n + 1)G^i -G^{m}_m y^i.  \label{C.2.2}
\end{equation}
Thus a Finsler space $F^n$ becpmes a Douglas space of the second kind if and only if (%
\ref{C.2.2}) are homogenous polynomials in $(y^i)$ of degree two.

\begin{definition}
A Finsler space $F^n$ is said to be a Douglas space of second kind if $%
D^{im}_m = (n+1)G^i -G^{im}_m y^i$ is a homogenous polynomial in ($y^i$) of
degree two.
\end{definition}

On the other hand, a Finsler space with an $(\alpha,\beta)$-metric is said
to be a $Douglas\ space\ of\ second\ kind$ if and only if 
\begin{equation}
B^{im}_m = (n + 1)B^i - B^m_m y^i,
\end{equation}
are homogenous polynomials in $(y^i)$ of degree two, where $B^m_m$ is given
by [\ref{CD-13}].

Furthermore, differentiating the above with respect to $y^h$, $y^j$ and $y^k$%
, we get 
\begin{equation}
B^{im}_{hjkm} = B^i_{hjk} = 0.
\end{equation}

\begin{definition}
A Finsler space $F^n$ with $(\alpha, \beta)$-metric is said to be a Douglas
space of second kind, if $B^{im}_m = (n + 1)B^i - B^m_m y^i$ is a homogenous
polynomial in $(y^i)$ of degree two.
\end{definition}

\section{Douglas Space of second kind with $(\protect\alpha, \protect\beta)$%
-metric}

In this section, we deal with the condition for a Finsler space with an $%
(\alpha, \beta)$-metric to be a Douglas space of second kind.

Let $G^i(x,y)$ be the spray function of a Finsler space $F^n$ with an $(\alpha, \beta)$%
-metric. According to [\ref{CD-7}], $G^i(x,y)$ is written in the form\newline
\begin{equation*}
2G^i = \gamma^i_{00} + 2 B^i,
\end{equation*}
\begin{equation}
B^i = \frac{\alpha L_\beta}{L_\alpha}s^i_0 + C^* \Bigl[\frac{\beta L_\beta}{%
\alpha L}y^i - \frac{\alpha L_{\alpha\alpha}}{L_{\alpha}} \Bigl(\frac{y^i}{%
\alpha} - \frac{\alpha b^i}{\beta} \Bigr)\Bigr],  \label{C.3.1}
\end{equation}
where 
\begin{eqnarray}
&&C^* = \frac{\alpha \beta (r_{00}L_\alpha - 2 \alpha s_0L_{\beta})}{2
(\beta^2 L_{\alpha} + \alpha \gamma^2 L_{\alpha \alpha})},  \notag \\
&& \gamma^2 = b^2 \alpha^2 - \beta^2.  \label{C.3.2}
\end{eqnarray}
Since $\gamma^i_{00} = \gamma^i_{jk}(x)y^jy^k$ is hp(2), equation (\ref%
{C.3.1}) yields 
\begin{equation}
B^{ij} = \frac{\alpha L_{\beta}}{L_{\alpha}}(s^i_0 y^j - s^j_0 y^i) + \frac{%
\alpha^2 L_{\alpha \alpha}}{\beta L_{\alpha}} C^*(b^iy^j - b^jy^i).
\label{C.3.3}
\end{equation}
By means of (\ref{C.2.1}) and (\ref{C.3.3}), we have the following lemma [%
\ref{CD-9}] :

\begin{lemma}
A Finsler space $F^n$ with an $(\alpha, \beta)$-metric becomes a Douglas space if
and only if $B^{ij} = B^i y^j - B^j y^i$ are hp(3).
\end{lemma}

Differentiating (\ref{C.3.3}) with respect to $y^h$, $y^k$, $y^p$ and $y^q$,
we have $D^{ij}_{hkpq} = 0$, which are equivalent of $D^{im}_{hkpm} = (n+1)
D^i_{hkp} = 0$. Thus, if a Finsler space $F^n$ satisfies the condition $%
D^{ij}_{hkpq} = 0$, we call it Douglas space. Further differentiating (\ref%
{C.3.3}) by $y^m$ and contracting m and j in the obtained equation, we
obtain 
\begin{eqnarray}
B^{im}_m &=& \frac{(n+1)\alpha L_{\beta} s^i_0}{L_{\alpha}} + \frac{\alpha
\{(n+1)\alpha^2 \Omega L_{\alpha \alpha}b^i + \beta \gamma^2 A y^i\}r_{00}}{%
2\Omega^2}  \notag \\
&&- \frac{\alpha^2 {\{(n+1)\alpha^2 \Omega L_{\beta}L_{\alpha \alpha}b^i +
By^i\}s_0}}{L_{\alpha} \Omega^2} - \frac{\alpha^3 L_{\alpha \alpha}y^i r_0}{%
\Omega},  \label{C.3.4}
\end{eqnarray}
where\newline
$\Omega$ = $(\beta^2 L_{\alpha} + \alpha \gamma^2 L_{\alpha \alpha}),\
provided\ that\ \Omega \neq 0$,\newline
$A$ = $\alpha L_{\alpha} L_{\alpha \alpha \alpha} + 3 L_{\alpha }L_{\alpha
\alpha} - 3 \alpha (L_{\alpha \alpha})^2$, 
\begin{equation}
B = \alpha \beta \gamma^2 L_{\alpha }L_{\beta}L_{\alpha \alpha \alpha} +
\beta \{(3\gamma^2 - \beta^2)L_{\alpha} -4\alpha \gamma^2 L_{\alpha
\alpha}\}L_{\beta}L_{\alpha \alpha } +\Omega L L_{\alpha \alpha}.
\label{C.3.5}
\end{equation}
We use the following result [\ref{CD-5}]:

\begin{theorem}
The necessary and sufficient condition for a Finsler space $F^n$ with an $%
(\alpha, \beta)$-metric to be a Douglas space of second kind is that, $%
B^{im}_m$ are homogenous polynomials in $(y^m)$ of degree two, where $%
B^{im}_m$ is given by (\ref{C.3.4}) and (\ref{C.3.5}), provided that $\Omega
\neq 0$.
\end{theorem}

\section{Conformal change of Douglas space of second kind with $(\protect%
\alpha, \protect\beta)$-metric.}

In the present section, we find the condition on conformal change, so that a
Douglas space of second kind with $(\alpha, \beta)$-metric is conformally transformed to a Douglas space
of second kind.

Let $F^n$= $(M^n, L)$ and \={F}$^n$ = ($M^n$, \={L}) be two Finsler spaces on
the same underlying manifold $M^n$. If we have a function $\sigma(x)$ in
each coordinate neighbourhoods of $M^n$ such that \={L}(x,y) = $e^\sigma
L(x,y)$, then $F^n$ is called conformal to \={F}$^n$ and the change L $%
\rightarrow$\={L} of metric is called a conformal change.

A conformal change of $(\alpha, \beta)$-metric is given as $(\alpha, \beta)$ 
$\rightarrow$ ($\bar{\alpha},\bar{\beta}$) , where $\bar{\alpha} =e^{\sigma}
\alpha$ and $\bar{\beta}= e^{\sigma}\beta$. Therefore, we have 
\begin{eqnarray}
\bar{a}_{ij} &=& e^{2\sigma} a_{ij}, \ \ \bar{b_i} = e^\sigma b_i \\
\bar{a}^{ij} &=& e^{-2\sigma} a^{ij}, \ \ \bar{b^i} = e^{-\sigma}b^i
\label{C.4.1}
\end{eqnarray}
and $b^2$ = $a^{ij}b_i b_j$ = $\bar{a}^{ij} \bar{b}_i\bar{b}_j$. Thus, we
state the following theorem for further use:\newline

\begin{theorem}
A Finsler space with $(\alpha, \beta)$-metric with the length b of $b_i$
with respect to the Riemannian metric $\alpha$ is invariant under any
conformal change of $(\alpha, \beta)$-metric.
\end{theorem}

From (\ref{C.4.1}), it follows that the conformal change of Christoffel
symbols is given by [\ref{CD-4}]: 
\begin{equation}
\bar{\gamma}_{jk}^i = \gamma_{jk}^i + \delta^i_j \sigma_k + \delta^i_k
\sigma_j - \sigma^i a_{jk},  \label{C.4.2}
\end{equation}
where $\sigma_j = \partial_j \sigma$ and $\sigma^i = a^{ij}\sigma_j$.

From (\ref{C.4.1}) and (\ref{C.4.2}), we have the following identities:%
\newline

\begin{equation*}
\bar{\triangledown}_j \bar{b_i} = e^{\sigma} ( \triangledown_j b_i + \rho
a_{ij} -\sigma_i b_j ),
\end{equation*}
\begin{equation*}
\bar{r}_{ij} = e^{\sigma} [r_{ij} + \rho a_{ij} - \frac{1}{2} (b_i \sigma_j
+ b_j \sigma_i)], 
\end{equation*}
\begin{equation*}
\bar{s}_{ij} = e^{\sigma} [s_{ij} + \frac{1}{2} (b_i \sigma_j -b_j
\sigma_i)],
\end{equation*}
\begin{equation*}
\bar{s}^i_j = e^{-\sigma}[s^i_j + \frac{1}{2} (b^i \sigma_j -b_j \sigma^i)], 
\end{equation*}
\begin{equation}
\bar{s}_j = s_j + \frac{1}{2}(b^2 \sigma_j - \rho b_j),  \label{C.4.3}
\end{equation}
where $\rho = \sigma_r b^r$.\newline
From (\ref{C.4.2}) and (\ref{C.4.3}), we can easily obtain the following: 
\begin{eqnarray}
\bar{\gamma}_{00}^i &=& \gamma_{00}^i + 2 \sigma_0 y^i - \alpha^2 \sigma_j,
\\
\bar{r}_{00} &=& e^{\sigma} (r_{00} + \rho \alpha^2 - \sigma_0 \beta), \\
\bar{s}^i_0 &=& e^{-\sigma}[s^i_0 + \frac{1}{2} (\sigma s_0 b^i - \beta
\sigma^i)], \\
\bar{s_0} &=& s_0 + \frac{1}{2} (\sigma_0 b^i - \rho \beta).  \label{C.4.4}
\end{eqnarray}
Next, we find the conformal change of $B^{ij}$ given in (\ref{C.3.3}). We
have $\bar{L}(\alpha, \beta)$ = $e^{\sigma} L(\alpha, \beta)$, and 
\begin{equation}
\bar{L}_{\bar{\alpha}} = L_{\alpha},\ \ \bar{L}_{\bar{\alpha}\bar{\alpha}}=
e^{-\sigma}L_{\alpha \alpha},\ \ \bar{L}_{\bar{\beta}} = L_{\beta},\ \ \bar{%
\gamma}^2= e^{2\sigma}\gamma^2.  \label{C.4.5}
\end{equation}
By using (\ref{C.3.2}), (\ref{C.4.4}), (\ref{C.4.5}) and theorem (3.1), we
obtain 
\begin{equation}
\bar{C}^* = e^{\sigma} (C^* + D^*),
\end{equation}
where 
\begin{equation}
D^* = \frac{\alpha\beta [(\beta \alpha^2 - \sigma_0 \beta)L_{\alpha} -
\alpha (b^2 \sigma_0 - \rho \beta)L_{\beta}]}{2(\beta^2 L_\alpha + \alpha
\gamma^2 L_{\alpha \alpha})}.
\end{equation}
Hence, under the conformal change $B^{ij}$ can be written as: 
\begin{eqnarray*}
\bar{B}^{ij} &=& \frac{\alpha L_{\beta}}{L_{\alpha}}(s^i_0y^j - s^j_0y^i) + 
\frac{\alpha^2 L_{\alpha \alpha}}{\beta L_{\alpha} }C^* (b^iy^j-b^jy^i) \\
&&+ \Bigl(\frac{\alpha \sigma_0 L_{\beta}}{L_{\alpha} }+ \frac{\alpha^2
L_{\alpha \alpha}}{\beta L_{\alpha}} D^*\Bigr) (b^iy^j-b^jy^i) - \frac{%
\alpha \beta L_{\beta}}{2 L_{\alpha}} (\sigma^i y^j - \sigma^j y^i), \\
&=& B^{ij} + C^{ij},
\end{eqnarray*}
where 
\begin{equation*}
C^{ij} = \Bigl(\frac{\alpha \sigma_0 L_{\beta}}{L_{\alpha} }+ \frac{\alpha^2
L_{\alpha \alpha}}{\beta L_{\alpha}} D^*\Bigr) (b^iy^j-b^jy^i) - \frac{%
\alpha \beta L_{\beta}}{2 L_{\alpha}} (\sigma^i y^j - \sigma^j y^i).
\end{equation*}
From (\ref{C.3.5}), we have 
\begin{equation}
\bar{\Omega} = e^{2\sigma}\Omega, \ \ \bar{A} = e^{-\sigma}A,\ \ \bar{B}=
e^{2\sigma}B.
\end{equation}
Now, we apply conformation transformation to $B^{im}_m$, and obtain 
\begin{equation}
\bar{B}^{im}_m = B^{im}_m + K^{im}_m,
\end{equation}
where 
\begin{eqnarray}
2 K^{im}_m &=& \frac{(n + 1) \alpha L_{\beta}}{L_\alpha}(\sigma_0 b^i -
\beta \sigma^i) + \alpha \Bigl\{ \frac{(n+1) \alpha^2 \Omega L_{\alpha
\alpha}b^i + \beta \gamma^2 A y^i}{\Omega^2}\Bigr\}(\rho \alpha^2 - \sigma_0
\beta)  \notag \\
&& - \Bigl[\frac{\alpha^2\{(n+1) \alpha^2 \Omega L_{\beta} L_{\alpha \alpha}
b^i + By^i\}}{L_{\alpha} \Omega^2} - \frac{\alpha^3 L_{\alpha \alpha} y^i}{%
\Omega} \Bigr] (b^2\sigma_0 - \rho \beta).  \label{C.4.7}
\end{eqnarray}
Thus, we have the following result:

\begin{theorem}
The necessary and sufficient condition for a conformal change of Douglas
space of the second kind to be a Douglas space of second kind, is that $%
K^{im}_m (x)$ are homogenous polynomial in ($y^i$) of degree two.
\end{theorem}

\section{Conformal change of Douglas space of second kind with special $(%
\protect\alpha, \protect\beta)$-metric $L = \protect\alpha +\protect\epsilon 
\protect\beta + k \frac{\protect\beta^2}{\protect\alpha}$}

Let us consider a Finsler space with special $(\alpha, \beta)$-metric  
\begin{eqnarray}
L = \alpha +\epsilon \beta + k \frac{\beta^2}{\alpha},
\end{eqnarray}
where $\epsilon$ and $k$ are constants.\newline
Then, from (29), we can easily find 
\begin{eqnarray}
L_{\alpha} &=& 1- \frac{k\beta^2}{\alpha^2},  \notag \\
L_\beta &=& \epsilon +\frac{ 2k\beta}{\alpha}, \\
L_{\alpha \alpha} &=& 2 k \frac{\beta^2}{\alpha^3},  \notag \\
L_{\alpha \alpha \alpha} &=& \frac{-6k\beta^2}{\alpha^4}.  \notag
\end{eqnarray}
Hence, from (\ref{C.3.5}), we have 
\begin{eqnarray}
\Omega &=& \frac{-3k\beta^4 + (1 + 2kb^2)\alpha^2 \beta^2}{\alpha^2},  \notag
\\
A &=& \frac{-12k^2 \beta^4}{\alpha^5}, \\
B &=& \frac{2k}{\alpha^6} \Bigl\{ (1 +2kb^2)\alpha^4 \beta^4 -6\epsilon k
b^2 \alpha^3 \beta^5 -2k (2 + 7kb^2)\alpha^2 \beta^6  \notag \\
&& + 6 \epsilon k \alpha \beta^7 + 15k^2 \beta^8\Bigr\}.  \notag
\end{eqnarray}
Thus, $K^{im}_m$ in (\ref{C.4.7}), reduces to 
\begin{equation}
2K^{im}_m = \frac{(n+1)(\epsilon \alpha^3 + 2k \alpha \beta)(\sigma_0 b^i -
\beta \sigma^i)}{(\alpha^2 -k\beta^2)} + p_1 + p_2 + p_3 + p_4,
\end{equation}
where 
\begin{eqnarray*}
p_1 &=& \frac{(n+1) b^i}{\{k\beta^4 - (1-2kb^2)\alpha^2 \beta^2\}^2}\Bigl\{%
2k\rho(1-2kb^2)\alpha^6 \beta^4 - 2k^2 \rho \alpha^4 \beta^6- 12k^2 \rho b^2
\alpha^2 \beta^7 y^i \\
&&- (1-2kb^2)2k \sigma_0 \alpha^4 \beta^5 + (\sigma_0 - 6 \rho y^i)2k^2
\alpha^2 \beta^7 + 12k^2 \sigma_0 b^2 \alpha^2 \beta^6 y^i -12k^2
\sigma_0\beta^8 y^i\Bigr\}, \\
p_2 &=& \frac{-12 k^2 \beta^5 \gamma^2(\rho \alpha^2 - \sigma_0 \beta)y^i}{%
\{(1+2kb^2)\alpha^2 \beta^2 - 3 k\beta^4\}^2}, \\
p_3 &=& \frac{2k(b^2 \sigma_0 \alpha^2 y^i - \rho \alpha^2 \beta y^i)}{%
(\alpha^2 -k\beta^2)\{(1+2kb^2)\alpha^2 \beta^2- 3k\beta^4\}^2}\Bigl\{ (1
+2kb^2)\alpha^4 \beta^4 -6\epsilon k b^2 \alpha^3 \beta^5 \\
&&-2k (2 + 7kb^2)\alpha^2 \beta^6 + 6 \epsilon k \alpha \beta^7 + 15k^2
\beta^8\Bigr\}, \\
p_4 &=& \frac{2kb^2 \alpha^2 y^i \sigma_0 - 2k\rho \alpha^2 \beta y^i}{%
(1+2kb^2)\alpha^2 - 3k \beta^2}.
\end{eqnarray*}
which shows that $K^{im}_m$ is a homogenous polynomial in $(y^i)$ of degree
two. Hence, we have the following theorem:\newline

\begin{theorem}
\label{Theorem 1-CD} A Douglas space of second kind with special $ \left( \alpha, \beta\right)  $-metric
$L = \alpha +\epsilon \beta + k \frac{\beta^2}{\alpha}$, where $%
\epsilon$ and k are constants, is conformally transformed to a Douglas space
of second kind.
\end{theorem}

From theorem \ref{Theorem 1-CD}, we can easily prove that a Douglas space of
second kind for certain Finsler spaces with $(\alpha, \beta )$-metric is
conformally transformed to a Douglas space of second kind. We have the
following cases.

\begin{enumerate}
\item[\textbf{Case (i).}] If $\epsilon$ = 1 and $k$ = 0, the special $ \left( \alpha, \beta\right)  $-metric reduces to $L = \alpha
+ \beta$, which is the well known Randers metric. In this case, $2K^{im}_m$ reduces to 
\begin{equation}
2K^{im}_m = (n+1)\alpha(\sigma_0 b^i - \beta \sigma^i),
\end{equation}
which shows that $K^{im}_m$ is a homogenous polynomial in $(y^i)$ of degree
two.\newline
Hence, we have the following:

\begin{corollary}
A Douglas space of second kind with Randers metric $L = \alpha +
\beta$, is conformally transformed to a Douglas space of second kind.
\end{corollary}

\item[\textbf{Case (ii).}] If $\epsilon$ = 0 and $k$ = 1, the special $ \left( \alpha, \beta\right) $-metric reduces to $L =
\alpha + \frac{\beta^2}{\alpha}$. In this case, $2k^{im}_m$ reduces to 
\begin{equation}
2K^{im}_m = \frac{(n+1)( 2 \alpha \beta)(\sigma_0 b^i - \beta \sigma^i)}{%
(\alpha^2 -\beta^2)} + q_1 + q_2 + q_3 + q_4,
\end{equation}
where 
\begin{eqnarray*}
q_1 &=& \frac{(n+1) b^i}{\{\beta^4 - (1-2b^2)\alpha^2 \beta^2\}^2}\Bigl\{%
2\rho(1-2b^2)\alpha^6 \beta^4 - 2 \rho \alpha^4 \beta^6- 12 \rho b^2
\alpha^2 \beta^7 y^i \\
&&- (1-2b^2)2 \sigma_0 \alpha^4 \beta^5 + (\sigma_0 - 6 \rho y^i)2 \alpha^2
\beta^7 + 12 \sigma_0 b^2 \alpha^2 \beta^6 y^i -12 \sigma_0\beta^8 y^i\Bigr\}%
, \\
q_2 &=& \frac{-12 \beta^5 \gamma^2(\rho \alpha^2 - \sigma_0 \beta)y^i}{%
\{(1+2b^2)\alpha^2 \beta^2 - 3 \beta^4\}^2}, \\
q_3 &=& \frac{2(b^2 \sigma_0 \alpha^2 y^i - \rho \alpha^2 \beta y^i)}{%
(\alpha^2 -\beta^2)\{(1+2b^2)\alpha^2 \beta^2- 3\beta^4\}^2}\Bigl\{ (1
+2b^2)\alpha^4 \beta^4 -2 (2 + 7b^2)\alpha^2 \beta^6 \\
&& + 15 \beta^8\Bigr\}, \\
q_4 &=& \frac{2b^2 \alpha^2 y^i \sigma_0 - 2\rho \alpha^2 \beta y^i}{%
(1+2b^2)\alpha^2 - 3\beta^2}.
\end{eqnarray*}
This shows that $K^{im}_m$ is a homogenous polynomial in $(y^i)$ of degree
two.\newline
Hence, we have the following:

\begin{corollary}
A Douglas space of second kind with special $(\alpha, \beta)$-metric $L = \alpha + \frac{\beta^2}{\alpha}$ is conformally transformed to a
Douglas space of second kind.
\end{corollary}

\item[\textbf{Case (iii).}] If $\epsilon$ = 1 and $k$ = 1, the special $ \left( \alpha, \beta\right)  $-metric reduces to the form $L =
\alpha+ \beta + \frac{\beta^2}{\alpha}$, which is the first approximate Matsumoto metric. In this case, $2K^{im}_m$ reduces to 
\begin{equation}
2K^{im}_m = \frac{(n+1)( \alpha^3 + 2 \alpha \beta)(\sigma_0 b^i - \beta
\sigma^i)}{(\alpha^2 -\beta^2)} + r_1 + r_2 + r_3 + r_4,
\end{equation}
where 
\begin{eqnarray*}
r_1 &=& \frac{(n+1) b^i}{\{\beta^4 - (1-2b^2)\alpha^2 \beta^2\}^2}\Bigl\{%
2\rho(1-2b^2)\alpha^6 \beta^4 - 2 \rho \alpha^4 \beta^6- 12 \rho b^2
\alpha^2 \beta^7 y^i \\
&&- (1-2b^2)2 \sigma_0 \alpha^4 \beta^5 + (\sigma_0 - 6 \rho y^i)2 \alpha^2
\beta^7 + 12 \sigma_0 b^2 \alpha^2 \beta^6 y^i -12 \sigma_0\beta^8 y^i\Bigr\}%
, \\
r_2 &=& \frac{-12 \beta^5 \gamma^2(\rho \alpha^2 - \sigma_0 \beta)y^i}{%
\{(1+2b^2)\alpha^2 \beta^2 - 3 \beta^4\}^2}, \\
r_3 &=& \frac{2(b^2 \sigma_0 \alpha^2 y^i - \rho \alpha^2 \beta y^i)}{%
(\alpha^2 -\beta^2)\{(1+2b^2)\alpha^2 \beta^2- 3\beta^4\}^2}\Bigl\{ (1
+2b^2)\alpha^4 \beta^4 -6 b^2 \alpha^3 \beta^5 \\
&&-2 (2 + 7b^2)\alpha^2 \beta^6 + 6 \alpha \beta^7 + 15\beta^8\Bigr\}, \\
r_4 &=& \frac{2b^2 \alpha^2 y^i \sigma_0 - 2\rho \alpha^2 \beta y^i}{%
(1+2b^2)\alpha^2 - 3\beta^2}.
\end{eqnarray*}
This shows that $K^{im}_m$ is a homogenous polynomial in $(y^i)$ of degree
two.\newline
Hence, we have the following:

\begin{corollary}
A Douglas space of second kind with first approximate Matsumoto metric $L
= \alpha+ \beta + \frac{\beta^2}{\alpha}$ is conformally transformed to a
Douglas space of second kind.
\end{corollary}

\item[\textbf{Case (iv).}] If $\epsilon$ = 2 and $k$ = 1, the special $ \left( \alpha, \beta\right)  $-metric reduces to the form $L = \frac{%
(\alpha+ \beta)^2}{\alpha}$, which is known as square metric. In this case, $%
2K^{im}_m$ reduces to 
\begin{equation}
2K^{im}_m = \frac{(n+1)(2 \alpha^3 + 2 \alpha \beta)(\sigma_0 b^i - \beta
\sigma^i)}{(\alpha^2 -\beta^2)} + s_1 + s_2 + s_3 + s_4,
\end{equation}
where 
\begin{eqnarray*}
s_1 &=& \frac{(n+1) b^i}{\{\beta^4 - (1-2b^2)\alpha^2 \beta^2\}^2}\Bigl\{%
2\rho(1-2b^2)\alpha^6 \beta^4 - 2 \rho \alpha^4 \beta^6- 12 \rho b^2
\alpha^2 \beta^7 y^i \\
&&- (1-2b^2)2 \sigma_0 \alpha^4 \beta^5 + (\sigma_0 - 6 \rho y^i)2 \alpha^2
\beta^7 + 12 \sigma_0 b^2 \alpha^2 \beta^6 y^i -12 \sigma_0\beta^8 y^i\Bigr\}%
, \\
s_2 &=& \frac{-12 \beta^5 \gamma^2(\rho \alpha^2 - \sigma_0 \beta)y^i}{%
\{(1+2b^2)\alpha^2 \beta^2 - 3 \beta^4\}^2}, \\
s_3 &=& \frac{2(b^2 \sigma_0 \alpha^2 y^i - \rho \alpha^2 \beta y^i)}{%
(\alpha^2 -\beta^2)\{(1+2b^2)\alpha^2 \beta^2- 3\beta^4\}^2}\Bigl\{ (1
+2b^2)\alpha^4 \beta^4 -12 b^2 \alpha^3 \beta^5 \\
&& -2 (2 + 7b^2)\alpha^2 \beta^6 +12 \alpha \beta^7 + 15 \beta^8\Bigr\}, \\
s_4 &=& \frac{2b^2 \alpha^2 y^i \sigma_0 - 2\rho \alpha^2 \beta y^i}{%
(1+2b^2)\alpha^2 - 3\beta^2},
\end{eqnarray*}
which shows that $K^{im}_m$ is a homogenous polynomial in $(y^i)$ of degree
two.\newline
Hence, we have the following:

\begin{corollary}
A Douglas space of second kind with square metric $L = \frac{(\alpha+
\beta)^2}{\alpha}$ is conformally transformed to a Douglas space of second
kind.
\end{corollary}
\end{enumerate}

\end{document}